\documentclass{amsart}
\usepackage{a4wide}

\usepackage{times}
\usepackage{amscd}
\usepackage{amssymb,amsmath}  
\usepackage{latexsym} 
\DeclareFontEncoding{OT2}{}{} 
\usepackage[english]{babel}

\newcommand{\Q}{\mathbf{Q}}
\newcommand{\A}{\mathbf{A}}
\newcommand{\Z}{\mathbf{Z}}

\newcommand{\C}{\mathbf{C}}

\newcommand{\Ps}{\mathbf{P}}

\newcommand{\ra}{\rightarrow}

\renewcommand{\bar}[1]{\overline{#1}}
\renewcommand{\phi}{\varphi}

    \newtheorem{Lem}{Lemma}[section]
    \newtheorem{Prop}[Lem]{Proposition}
    \newtheorem{Thm}[Lem]{Theorem}
    \newtheorem{Cor}[Lem]{Corollary}

   \theoremstyle{definition}
    \newtheorem{Def}[Lem]{Definition}

    \newtheorem{Rem}[Lem]{Remark}
    
  \newtheorem{Results}[Lem]{Results}
    \DeclareMathOperator{\Pic}{Pic}
    \DeclareMathOperator{\Gal}{Gal}
\DeclareMathOperator{\rank}{rank}
\DeclareMathOperator{\id}{id}

\DeclareMathOperator{\Jac}{Jac} 
\DeclareMathOperator{\Aut}{Aut}

\DeclareMathOperator{\res}{res}

\DeclareMathOperator{\End}{End}

\begin{document}
\title[Kuwata's surfaces]{Explicit sections on Kuwata's elliptic surfaces}
\author{Remke Kloosterman}
\address{Department of Mathematics and Computer Science, University of Groningen, PO Box 800, 9700 AV  Groningen, The Netherlands}
\email{r.n.kloosterman@math.rug.nl}
\thanks{The author would like to thank Jaap Top and Marius van der Put for many valuable discussion on this paper. This paper is based on a chapter in the author's PhD thesis \cite[Chapter 3]{Proefschrift}.
}
\begin{abstract}
We give explicit generators for subgroups of finite index of the Mordell-Weil groups of several families of elliptic surfaces introduced by Masato Kuwata.
\end{abstract}

\keywords{Elliptic surfaces}
\date{\today}
\maketitle

\section{Introduction}
In this article we assume that $K$ is a field of characteristic 0.

We introduce elliptic surfaces $\pi_n: X_n \ra \Ps^1$ defined over $K$, for every integer $n\geq 1$. In the special case $n\equiv 0 \bmod 2$,  the function field extension $K(X_n)/K(\Ps^1)$ is $K(x,z,t)/K(t)$ where $x^3+ax+b-t^n(z^3+cz+d)=0$, for some $a,b,c,d\in K$. For every $n\geq 1$ the surface $X_n$ is birational to a base-change $X_1 \times_{\Ps^1} \Ps^1$, where $\Ps^1\ra \Ps^1$ is an $n$-cyclic cover.

Kuwata \cite{Kuw} computed the rank of the Mordell-Weil groups $MW(\pi_i):=MW_{\bar{K}}(\pi_i)$ of $\pi_i$, for $i=1,\dots, 6$, and he gave a strategy for computing a rank 12 subgroup of the rank 16 group $MW(\pi_6)$. This strategy for finding generators of $MW(\pi_6)$ is already mentioned in an email correspondence between Jasper Scholten and Masato Kuwata.  We  describe and extend their ideas in order to describe  generators of the Mordell-Weil groups $MW(\pi_i)$ for $i=2,3,4$ and 6. Most of the computations we use, can be found in the Maple worksheet \cite{Worksheet}. 
We indicate how one can find in special cases generators for the Mordell-Weil group of $\pi_5$. Finally, in the case $K=\Q$ we discuss how large the part of the Mordell-Weil group consisting of $\Q$-rational sections can be.

For $i\leq 6$ the $\pi_i:X_i\ra \Ps^1$ are elliptic $K3$ surfaces. There are very few examples of such surfaces where 
generators for the Mordell-Weil group have been found.

\section{Notation and results}\label{notres}
Fix two elliptic curves $E$ and $F$, with $j(E)\neq j(F)$. Let $\iota: E\times F\ra E\times
F$ be the automorphism sending $(P,Q)\ra (-P,-Q)$. The minimal
desingularization of $(E\times F)/\langle \iota\rangle$ is a
$K3$ surface which we denote by $Y$. It is called the
Kummer surface of $E\times F$.

The surface $Y$ possesses several elliptic fibrations. For a generic
choice of $E$ and $F$, all fibrations on $Y$ are classified by
Oguiso (\cite{Ogui}). Following Kuwata \cite{Kuw}, we concentrate on a particular fibration $\psi:
Y \ra \Ps^1$ having two fibers of type $IV^*$ and 8 other irreducible singular fibers. Assume that the fibers of type $IV^*$ are over $0$ and $\infty$. 
Let $\pi_6:X_6 \ra \Ps^1$ be the cyclic degree 3 base-change of $\psi$ ramified over $0$ and $\infty$. Assume that $E$ is given by $y^2=x^3+ax+b$ and $F$ is given by $y^2=x^3+cx+d$. Set $B(t)=\Delta(F) t+864bd+\Delta(E)/t$, with $\Delta(E)$ the discriminant of $E$, i.e., $\Delta(E)=-16(4a^3+27b^2)$, and $\Delta(F)=-16(4c^3+27d^2)$. Then a Weierstrass equation for $\pi_6$ is
\[ y^2=x^3-48ac x+B(t^6). \]
Define $\pi_i:X_i\ra \Ps^1$ as the elliptic surface associated to the Weierstrass equation
\[ y^2=x^3-48ac x+B(t^i).\]
Clearly, interchanging the role of $E$ and $F$ corresponds to the automorphism $t\mapsto 1/t$ on $X_i$. One has that $j(\pi_2)=j(\psi)$, but $X_2$ and $Y$ are {\em not} isomorphic as fibered surfaces.

By construction, it is clear that $MW(\pi_n)$ can be regarded as a subgroup of $MW(\pi_{nm})$ for $m\geq 1$.

We recall the following result from \cite{Kuw}, most of which will be reproven in the course of Section~\ref{results}.
\begin{Thm}[{Kuwata, \cite[Theorem 4.1]{Kuw}}]\label{exactrank}  The surface $X_i$ is a $K3$ surface if and only if
  $i\leq 6$. Set 
\[ h=\left\{\begin{array}{cl} 0 &  \mbox{if $E$ and $F$ are not isogenous,}\\ 1&\mbox{if $E$ and $F$ are isogenous and $E$ does not admit complex multiplication,} \\ 2&\mbox{if $E$ and $F$ are isogenous and $E$ admits complex multiplication.}\end{array}\right.\]
Suppose $j(E)\neq j(F)$, then
\[ \rank MW(\pi_i)=\left\{\begin{array}{cl} h &\mbox{if }  i=1,\\4+h&\mbox{if } i=2,\\8+h&\mbox{if } i=3,\\12+h&\mbox{if } i=4,\\ 16+h & \mbox{if } i=5, \\ 16+h & \mbox{if } i=6. \end{array}\right.\]
\end{Thm}
In the case that $j(E)=j(F)$ then the ranks of $MW(\pi_i)$ tend to be lower (see \cite[Theorem 4.1]{Kuw}).
Kuwata in \cite{Kuw} restricts himself to the $X_i$ with $i\leq 6$. We observe the following corollary of his results.
\begin{Cor}\label{minrank}
The rank of $ MW(\pi_{60})$ is at least $ 40+h$.
\end{Cor} 

\begin{proof}It can be easily seen that the rank $40+h$ group\[ MW(\pi_1) \oplus \frac{MW(\pi_2)}{MW(\pi_1)}\oplus\frac{
  MW(\pi_4)}{MW(\pi_2)} \oplus \frac{MW(\pi_5)}{MW(\pi_1)}\oplus\frac{
  MW(\pi_6)
}{MW(\pi_2)}.\]
injects into $MW(\pi_{60})$. (The summands correspond to different eigenspaces for the induced action of $t \mapsto \zeta_{60} t$ on $MW(\pi_{60})$.)
\end{proof}
\begin{Rem} The highest known
  Mordell-Weil rank for an elliptic surface $\pi: X \ra \Ps^1$, such
  that $j(\pi)$ is non-constant is 56, a result due to Stiller \cite{Stil}. The geometric genus of our example is very low compared to for example Stiller's examples: one can easily
  show that $p_g(X_{60})=19$, while  Stiller's examples have $p_g+1$ divisible by 210.\end{Rem}

Our aim is to provide explicit equations for generators of $MW(\pi_i)$, for  several small values of $i$, in terms of the parameters $a,b,c,d$ of the two elliptic curves
\[ E: y^2=x^3+ax+b \mbox{ and } F:y^2=x^3+cx+d.\]
The main idea used in the sequel is the following. Let $\pi: S \ra \Ps^1$ be a $K3$ surface. Let $\sigma$ be an automorphism of finite order of $S$, such that we have a commutative diagram
\[ \begin{array}{lcl}
S& \stackrel{\sigma}{\ra} &S\\
\downarrow\pi &&\downarrow\pi\\
\Ps^1 &\stackrel{\tau}{\ra}& \Ps^1\end{array}\]
and the order of $\tau$ equals the order of $\sigma$. We obtain an elliptic surface $\psi:S':=\widetilde{X/\langle \sigma\rangle} \ra \Ps^1/\langle \tau\rangle \cong \Ps^1$. In the case that $S'$ is a rational surface, it is in principle possible to find explicit equations for generators of $MW(\psi)$. One can pull back sections of $\psi$ to sections of $\pi$, which establishes $MW(\psi)$ as a subgroup of $MW(\pi)$. In the case that $\sigma$ is of even order then there exists a second elliptic surface, $\psi'': S'' \ra \Ps^1$, such that $MW(\psi) \oplus MW(\psi'')$ modulo torsion injects into $MW(\pi)$ modulo torsion. (This process is called twisting and is discussed in more detail in Section~\ref{Prel}.)

In \cite[Section 5]{Kuw} it is indicated how one can find generators for a rank 12 subgroup of $MW(\pi_6)$ by using the above strategy for the involutions $(x,y,t)\mapsto (x,y, -t)$, and $(x,y,t)\mapsto (x,y,\alpha/t)$, for some $\alpha$ satisfying $\alpha^3=\Delta(F)/\Delta(E)$. Taking a third involution $(x,y,t) \mapsto (x,y,\alpha\zeta_3/t)$ yields a rank 16 subgroup of $MW(\pi_6)$, isomorphic to $MW(\pi_6)/MW(\pi_1)$. (This is discussed in more detail in the Subsections~\ref{ssecthree} and~\ref{ssecsix}.) In Subsection~\ref{ssecsix} we also describe the field of definition of sections generating  $MW(\pi_6)$ modulo $MW(\pi_1)$. It turns out that our  description corrects a mistake in Kuwata's description of the minimal field of definition.

\begin{Results}\label{resu}On the  Kuwata surfaces $X_1,X_2,\dots, X_6$ one has several automorphisms such that the obtained quotients are rational surfaces. This is used in Section~\ref{results} to give the following results:

\begin{itemize}
\item We give sufficient conditions on $E, F$ and $K$ to have a rank 1 or a rank 2 subgroup of $MW_K(\pi)$.
\item We give an indication how one can find explicit generators of $MW(\pi_2)/MW(\pi_1)$. The precise generators can be found in \cite{Worksheet}.
\item We give an indication how one obtains a degree 24 polynomial such that its zeroes determine a set of generators of $MW(\pi_3)/MW(\pi_1)$. The explicit polynomial can be found in \cite{Worksheet}.
\item We give explicit generators of a subgroup of finite index in $MW(\pi_4)/MW(\pi_2)$. (See Corollary~\ref{res4}.)
\item We give an indication how one obtains a degree 240 polynomial such that its zeroes determine a set of generators of $MW(\pi_5)/MW(\pi_1)$. 
We did not manage to write down this polynomial. The degree of this polynomial and the number of variables is too high to write it down explicitly. For several choices of $a,b,c,d$  it can be found in \cite{Worksheet}.
\item We give an algorithm for determining a set of generators of  $MW(\pi_6)/MW(\pi_3)$. The explicit generators are given in \cite{Worksheet}.
\end{itemize}\end{Results}

If $E$ and $F$ are not isogenous then $MW(\pi_1)=0$ (see Theorem~\ref{exactrank}). If $E$ and $F$ are isogenous then the degree of the $x$-coordinate of a generator of $MW(\pi_1)$ seems to depend on the degree of the isogeny between $E$ and $F$. This seems to give an obstruction for obtaining explicit formulas for all cases.

Before providing the explicit equations, we mention the following result, which can be proven without knowing explicitly the generators of $MW(\pi_i)$.

\begin{Prop}\label{qrank} Suppose that $E$ and $F$ are defined over $\Q$. One has that 
\[ \rank MW_{\Q}(\pi_i)\leq\left\{
\begin{array}{cl} 1 &\mbox{if } i=1,\\5&\mbox{if } i=2,\\7&\mbox{if } i=3,\\9&\mbox{if } i=4,\\5&\mbox{if } i=5,\\11&\mbox{if } i=6.\end{array}\right.\]
The  total contribution of $\pi_1, \pi_2,\pi_3,\pi_4,\pi_6$ to $MW_{\Q}(\pi_{12})$ is bounded by 15. The maximum total contribution of $\pi_1,\dots,\pi_6$ to $MW_{\Q}(\pi_{60})$ is bounded by 19.
\end{Prop}
\begin{proof}
If $\rank MW_{\Q}(\pi_1)$ equals  2 then  the Shioda-Tate formula~\ref{ST} implies that $NS(X)$ has rank 20 and $NS(X)$ is generated by divisors defined over $\Q$. It is well-known that this is impossible (\cite{no20}).

Suppose that $i\geq 2$. Let $S$ be the image of a section  of $\pi_i$, such that $S$ is not the strict transform of the  pull back of a section of $\pi_j$, for some $j$ dividing $i$,  and $j\neq i$.
It is easy to check that if we push forward $S$ to $X_1$ and then pull this divisor  back to $X_i$, we obtain a divisor $D$ consisting of $i$ geometrically irreducible components, and one of these components is $S$.

Set $M_j:=MW(\pi_j)$. Set $M:=M_i/\sum_{j>0,j|i,j\neq i} M_j$.
Take a minimal set of sections $S_m, m=1,\dots,\ell $, satisfying the following property: denote $S_{m,n}$ the components of the  pull back of the push forward of $S_m$, then the $S_{m,n}$, with $m=1,\dots, \ell$ and $n=1,\dots,i$, generate $M\otimes \Q$. 
Consider the $\ell n$-dimensional $\Q$-vector space $F$ of formal linear expression in the $S_{m,n}$. By definition of the $S_n$ it follows that the natural map $\Psi: F \ra M \otimes \Q$
 is surjective. Fix a generator $\sigma$ of the (birational) automorphism group of the rational map $X_i\ra X_1$. We can split $F\otimes \C$ and $M_j \otimes \C$ into eigenspaces for the action of $\sigma$. Set $E_{\zeta_i^k}\subset F \otimes \C$ the eigenspace for the eigenvalue $\zeta_i^k$. Set $V':=\oplus_{j \mid \gcd(i,j)=1} E_{\zeta_i^j}$. Since the $S_{i,j}$ form a minimal set, we have that that $\Psi|_{V'}$ is injective. It is easy to see that if $\gcd(i,k)\neq 1$ then the eigenspace $E_{\zeta_i^k}$ is contained in the kernel of $\Psi$, hence we obtain that $\Psi|_{V'}$ is an isomorphism.
This implies that $\dim V'=\varphi(i)\ell$, where $\varphi$ is the Euler $\varphi$-function.

Fix $m$ such that $1\leq m \leq  \ell$. Let $G$ be the subgroup of $MW(\pi_i)$ generated by the $S_{m,n}$ for $n=1,\dots i$. There is a faithful action of the Galois group $\Gal(\Q(\zeta_i)/\Q)$ on $G$.   
This implies that the sections defined over $\Q$ form  a rank at most 1 subgroup of $G$  (if $i$ is odd) or a rank at most two subgroup of $G$  (if $i$ is even). From this one obtains that
\[M_{\Q}:= MW_{\Q}(\pi_i)/\oplus_{j>0,j|i,j\neq i}MW_{\Q}(\pi_j)\]
has rank at most $\ell$ (when $i$ is odd) or $2\ell$ (when $i$ is even). This combined with Proposition~\ref{exactrank} provides the upper bounds for $MW_{\Q}(\pi_i)$ for $i=2,3,5,6$. The upper bound for $MW_{\Q}(\pi_4)$ follows from the fact that the above mentioned Galois action is non-trivial, hence for each $i$ there is an element $g$ in $G$ not fixed under the action of $\Gal(\Q(\zeta_i)/\Q)$ and $g$ is not mapped to $0$ in $M$. One easily obtains from this that $M_{\Q}$ has co-rank at least $\ell$ in $M$, which gives the case $i=4$.

The final assertions follow from a reasoning as in the proof of Corollary~\ref{minrank}.
\end{proof}

For some of the $\pi_i$ the bounds in Proposition~\ref{qrank} are sharp, for some of the others we do not know:
\begin{itemize}
\item  We have $\rank MW_{\Q}(\pi_1)=1$ when $E$ and $F$ are isogenous over $\Q$. (see \ref{ssecone}.)
\item We have that $\rank MW_{\Q}(\pi_2)\geq 4$ if and only if $E$ and $F$ have complete two-torsion over $\Q$, and, moreover, if $E$ and $F$ are isogenous then $\rank MW_{\Q}(\pi_2)=5$. (see \ref{ssectwo}.)
\item There exists $E$ and $F$ such that $\rank MW_{\Q}(\pi_4)\geq 8$. (see \ref{ssecfour}.)
\item There exists $E$ and $F$ such that $\rank MW_{\Q}(\pi_6)\geq 6$.

\item If one can find a rational solution of $uv^2(u^4-1)(v^4-1)^2=w^3$, satisfying
\begin{itemize}
\item $w(u-v)\neq 0$ and
\item if $p_1=4u^2/(u^2-1)$ and $p_2=4v^2/(v^2-1)$ then $p_1 \not \in \{p_2,1/p_2,1-p_2,1-1/p_2, p_2/(1-p_2),1/(1-p_2)\}$.
\end{itemize} then one can make examples with $\rank MW_{\Q}(\pi_{12})\geq 10$.
(see \ref{ssecsix}.)
\end{itemize}

\begin{Rem} The highest known rank over $\Q$ for an elliptic surface $\pi: X\ra\Ps^1$  over $\Q$ is  14 \cite{Kihara}. This example is a member of a family of elliptic surface introduced by Mestre \cite{Mestre}.\end{Rem}

The actual aim of this article is to produce explicit equations for the
generators of $MW(\pi_i)$ modulo $MW(\pi_1)$, for $i=2,\dots, 6$. 
For convenience, we always assume
that $E$ and $F$ have complete 2-torsion over the base field. It is not so hard to deduce from our results, these equations for the case that $E$ and $F$ have 2-torsion points defined over a larger field.

The organization of this paper is as follows. In Section~\ref{Prel} we recall some standard definitions. In Section~\ref{secfind} we indicate an algorithm for finding generators of $MW(\pi)$ in the case that of a rational elliptic surface $\pi:X\ra\Ps^1$  with an additive fiber. In Section~\ref{results} we prove the results stated in~\ref{resu}. To find generators for the $MW(\pi_i)$ we often rely on the results of Section~\ref{secfind}.

\section{Definitions}
\label{Prel}

\begin{Def}\label{defbas} An \emph{elliptic surface} is a triple
  $(\pi,X,C)$ with $X$ a surface, $C$ a curve, $\pi$ is a morphism
  $X\rightarrow C$, such that almost all fibers are irreducible genus
  1 curves and $X$ is relatively minimal, i.e., no fiber of $\pi$
  contains an irreducible rational curve $D$ with $D^2=-1$.

  We denote by $j(\pi): C \rightarrow \Ps^1$ the rational function such that $j(\pi)(P)$ equals the $j$-invariant of $\pi^{-1}(P)$, whenever $\pi^{-1}(P)$ is non-singular.

A \emph{Jacobian elliptic surface} is an elliptic surface together with a section $\sigma_0: C \rightarrow X$ to $\pi$.
  The set of sections of $\pi$ is an abelian group, with $\sigma_0$ as the identity element. Denote this group by $MW(\pi)$.

Let $NS(X)$ be the group of divisors on $X$ modulo algebraic equivalence, called the {\em N\'eron-Severi group} of $X$.
Let  $\rho(X)$ denote the rank of the
  N\'eron-Severi group of $X$. We call $\rho(X)$ the {\em Picard number}.
\end{Def}

\begin{Def} Let $X$ be a surface, let $C$ and $C_1$ be curves. Let $\varphi: X \rightarrow C$ and $f: C_1 \rightarrow C$ be two morphisms. Then we denote by $\widetilde{X \times_C C_1}$ the smooth, relatively minimal model of the ordinary fiber product of $X$ and $C_1$.
\end{Def}

\begin{Rem}\label{singtabel} If $P$ is a point on $C$, such that $\pi^{-1}(P)$ is
  singular then $j(\pi)(P)$ and $v_p(\Delta_p)$ behave as in Table~\ref{kod}.
For proofs of these facts see \cite[p. 150]{BPV}, \cite[Theorem IV.8.2]{Silv2} or \cite[Lecture 1]{MiES}.
\end{Rem}
\begin{table}
\[ \begin{array}{|c|c|c|}
\hline
\mbox{Kodaira type of fiber over }P & j(\pi)(P)  & \mbox{ number of components}\\
\hline
I_0^*                  & \neq \infty & 1\\
I_\nu \;(\nu>0) & \infty &\nu+1\\
I_\nu^* \;(\nu>0) & \infty &\nu+5 \\
II     & 0 & 1\\
IV &0& 3\\
IV^* & 0 & 7 \\
II^* & 0 &9\\
III &1728 & 2\\
III^*              & 1728 &8\\ \hline \end{array}\]
\caption{Classification of singular fibers}\label{kod}\end{table}

Recall the following theorem.
\begin{Thm}[{Shioda-Tate (\cite[Theorem 1.3 \& Corollary
    5.3]{Sd})}]\label{ST}  Let $\pi:X\rightarrow C$ be a Jacobian
  elliptic surface, such that $\pi$ has at least one singular fiber. Then the N\'eron-Severi
  group of $X$ is generated by the classes of $\sigma_0(C)$, a
  non-singular fiber, the components of the singular fibers not
  intersecting $\sigma_0(C)$, and the generators of the Mordell-Weil
  group. Moreover, let $S$ be the set of points $P$ such that $\pi^{-1}(P)$ is singular. Let $m(P)$ be the number of irreducible components of $\pi^{-1}(P)$, then
  \[ \rho(X) =2+ \sum_{P \in S} (m(P)-1)+\rank(MW(\pi)) \]
\end{Thm}

The following result will be used several times. It is a direct
consequence of the Shioda-Tate formula.
\begin{Thm}[{\cite[Theorem 10.3]{Sd}}]\label{ratTHM} Let $\pi : X\ra
  \Ps^1$ be a rational Jacobian elliptic surface, then the rank of the Mordell-Weil group is 8 minus the number of irreducible components of singular fibers not intersecting the identity component. \end{Thm}

\begin{Def} Suppose $\pi: X \ra C$ is an elliptic surface. Denote by $T(\pi)$ the  subgroup of the N\'eron-Severi group of $\Jac(\pi)$ generated by the classes of the fiber, $\sigma_0(C)$ and the components of the singular fibers not intersecting $\sigma_0(C)$. Let $\rho_{tr}(\pi):=\rank T(\pi)$. We call $T(\pi)$ the {\em trivial part} of the N\'eron-Severi group of $\Jac(\pi)$.
\end{Def}

Given a Jacobian elliptic surface  $\pi: X \ra C$ over a field $K$, we can associate an
 elliptic curve in $\Ps^2_{K(C)}$  corresponding to the generic fiber
 of $\pi$. This induces a bijection on isomorphism classes of Jacobian
 elliptic surfaces and elliptic curves over $K(C)$.

Two elliptic curves $E_1$ and $E_2$ are isomorphic over $K(C)$ if and
only if $j(E_1)=j(E_2)$ and the quotients of the minimal discriminants of $E_1/K(C)$ and $E_2/K(C)$ is a 12-th power (in $K(C)^*$).

Assume that $E_1$, $E_2$ are elliptic curves over $K(C)$ with
$j(E_1)=j(E_2)\neq 0,1728$. Then one shows easily that
$\Delta(E_1)/\Delta(E_2)$ equals $u^6$, with $u\in K(C)^*$. Hence $E_1$ and
$E_2$ are isomorphic over $K(C)(\sqrt{u})$. We call $E_2$ the twist
of $E_1$ by $u$, denoted by $E_1^{(u)}$. Actually, we are not interested in the function $u$, but in the places at which the valuation of $u$ is odd.

\begin{Def} Let $\pi:X \ra C$ be a Jacobian elliptic surface. Fix $2n$ points $P_i \in C(\C)$. Let $E/\C(C)$ be the Weierstrass model of the generic fiber of $\pi$.

A Jacobian elliptic surface $\pi': X' \ra C$ is called a \emph{(quadratic) twist} of $\pi$ by $(P_1,\ldots,P_n)$ if the Weierstrass model of the generic fiber of $\pi'$ is isomorphic to $E^{(f)}$, where $E^{(f)}$ denotes the quadratic twist of $E$ by $f$ in the above mentioned sense
and $f\in\C(C)$ is a function such that  $v_{P_i}(f) \equiv 1 \bmod 2$ and $v_Q(f) \equiv 0 \bmod 2$ for all $Q\not \in \{P_i\}$.
\end{Def}
The existence of a twist of $\pi$ by $(P_1,\ldots,P_{2n})$ follows
directly from the fact that $\Pic^0(C)$ is 2-divisible.

If  we fix $2n$ points $P_1,\ldots P_{2n}$ then there
exist precisely $2^{2g(C)}$ twists by $(P_i)_{i=1}^{2n}$.

If $P$ is one of the $2n$ distinguished points, then the fiber of $P$ changes in the following way (see \cite[V.4]{MiES}).
\[
I_\nu \leftrightarrow I^*_\nu \;(\nu \geq 0) \;\;\; \;\;
II \leftrightarrow    IV^* \;\;\;\;\;
III  \leftrightarrow  III^* \;\;\;\;\;
IV  \leftrightarrow  II^*
\]

Let $\pi: X\ra  C$ be a Jacobian elliptic surface, $P_1,\dots P_{2n}\in C$
points. Let $\tilde{\pi}: \tilde{X}\ra C$ be a twist by the $P_i$. Let
$\varphi: C_1\ra C$ be a double cover ramified at the $P_i$, such that the
minimal models of base-changing $\varphi$ and $\tilde{\varphi}$  by
$\pi$ are isomorphic. Denote this model by $\pi_1:X_1\ra C_1$.

Recall that
\begin{eqnarray} \rank(MW(\pi_1)) = \rank(MW(\pi))+\rank(MW(\tilde{\pi})).\label{twteqn}\end{eqnarray} 
Moreover, the singular fibers change as follows
\[ \begin{array}{|l|cccc|}
\hline
\mbox{Fiber of } \pi  \mbox{ at } P_i                 & I_\nu \mbox{ or
  }I_\nu^*& II \mbox{ or } IV^*&
III  \mbox{ or } III^*& IV \mbox{ or } II^* \\

\hline
\mbox{Fiber of } \pi_1 \mbox{ at } \varphi^{-1}(P_i) & I_{2\nu}& IV &
I_0^* & IV^*\\\hline\end{array} \]

\section{Finding sections on a rational elliptic surface with an additive fiber}\label{secfind} 
Fix a field $K$ of characteristic 0 and a rational elliptic surface $\pi: X \ra \Ps^1$ over $K$. Such a surface can be
 represented by a Weierstrass equation
\begin{equation}\label{Weieqn} y^2=x^3+(\alpha_4t^4+\alpha_3t^3+\alpha_2t^2+\alpha_1t+\alpha_0)x+\beta_6t^6+\beta_5t^5+\beta_4t^4+\beta_3t^3+\beta_2t^2+\beta_1t+\beta_0,\end{equation}
with $\alpha_i,\beta_j\in K$. Assume that over $t=\infty$ the fiber is singular and of additive type. This happens if and only if $\alpha_4=\beta_6=0$.

From \cite[Theorem 2.5]{OS} we know that $MW(\pi)$ is generated by sections of the form
\begin{equation}\label{seceqn}  x=b_2t^2+b_1t+b_0, y=c_3t^3+c_2t^2+c_1t+c_0.
\end{equation}

Substituting (\ref{seceqn}) in (\ref{Weieqn}) and using $\alpha_4=\beta_6=0$ yields that
\[ c_3^2=b_2^3.\]
Set  $c_3=p_1^3,b_2=p_1^2$. First we search for solutions with $p_1\neq 0$. The equation for the coefficient of $t^5$ is an equation of the from $p_1^2c_2+f$ with $f$ a polynomial in $p_1, b_1, \alpha_i, \beta_j$. Hence we can express $c_2$ in terms of the $p_1,b_1, \alpha_i, \beta_j$. Similarly, we can express $c_1$ and $c_0$ in terms of $p_1,b_k, \alpha_i,\beta_j$. One easily shows that if this procedure  fails, then $p_1=0$.

Unfortunately, the three remaining equations $F_m=0$ in $b_1,b_0,p_1, \alpha_i,\beta_j$ are not linear in  any of the variables, but the degree of the $F_i$ in $b_1,b_0,p_1$ is sufficiently low to compute the resultant
\[ R(p_1):=\res_{b_1}(\res_{b_0}(F_1,F_3),\res_{b_0}(F_2,F_3))\]
in concrete examples (i.e., after substituting constants for  the $\alpha_i, \beta_j$). 
Calculating $R(p_1)$ and finding zeroes of it gives all possibilities for $p_1$. Substituting such a value  of $p_1$ in $\res(F_1,F_3,b_0)$ and in $\res(F_2,F_3,b_0)$ yields two polynomials in $b_1$. Calculating the g.c.d. of these polynomials  gives a list of possible values for $b_1$. Then substituting all possibilities for $(p_1,b_0)$ in $F_1,F_2,F_3$ gives all possibilities for $p_1,b_0,b_1$. 

Consider the case $p_1=0$, hence  $c_3=b_2=0$. The coefficient of $t^5$ is zero if and only if  $\alpha_5=0$. If this is the case then the coefficient of $t^4$ is of the form $-c_2^2+\beta_3b_1+\alpha_4$. 

If $\beta_3\neq 0$, then one eliminates  $b_1,b_0,c_0$ as is done above, yielding  two polynomial equation in two unknowns. In this case it suffices to compute only one resultant.

If $\beta_3=0$ then we can substitute $c_2=\pm \sqrt{\alpha_4}$. The coefficient of $t^3$ is of the form $-2\sqrt{\alpha_4}c_1+f$, with $f$ a polynomial in $c_0,\alpha_i,\beta_j,\sqrt{\alpha_4}$. We can solve $c_1$ and obtain  $-2\sqrt{\alpha_4}c_0+g$ as a coefficient for $t^3$, where $g$ is a polynomial in $c_0,\alpha_i,\beta_j,\sqrt{\alpha_4}$.

This fails when $\alpha_4=0$. In that case $c_2=0$, and one has four polynomial relations in $b_0,c_0,c_1,b_1$. We can eliminate $c_1$ as above, yielding three polynomial relations in three unknowns which can be solved as above.
In this way we find all possible sections of the form (\ref{seceqn}).

Top \cite[Section 5]{Top3des} discusses the following elliptic surface:
\begin{Thm} Let $\pi: X \ra \Ps^1$ be the elliptic surface
\[ y^2=x^3+108(27t^4-74t^3+84t^2-48t+12).\]
Then $MW_{\Q}(\pi)$ has rank 3 and is generated by sections $\sigma_i$
 with $x$-coordinates
\[x(\sigma_1)=6t,\; x(\sigma_2)=6t-8, \; x(\sigma_3)=-12t+9.\]
and fix a primitive cube root of unity $\omega$. Let $\tau_i$ be obtained from $\sigma_i$ by multiplying the $x$-coordinate with $\omega$. Then the $\sigma_i$ and $\tau_i$ generate a subgroup of finite index of $MW(\pi)$.
\end{Thm}

 Top found the sections $\pm \sigma_1,\pm \sigma_2, \pm \tau_1$ and $\pm \tau_2$; an explicit description of the third independent section $\sigma_3$ seems not to be present in the literature.

\begin{proof}
The elliptic surface $\pi: X \ra \Ps^1$
 has 4 fibers of type $II$ and at $t=\infty$ a fiber of type $IV$. 
From this it follows that $MW(\pi)$ has rank 6.
Top \cite[Section 5]{Top3des} observes that the Mordell-Weil group is generated by a subset of the 18 sections of the form $(\omega x_i,\pm y_i), i=1,2,3$, with $\omega^3=1$ and $x_i$ is a polynomial of degree 1.

 The form of the generators  imply that, in terms of the above discussion, we are looking for sections with $p_1=0$. After eliminating all the $c_i$ and $b_j$, except for $b_1$, we obtain a polynomial of degree 27. It can be factored as the product of $(b_1-6)^2(b_1+12)$ times 12 polynomials of degree 2, where all the degree 2 factors have discriminant $-3$.

This implies that $x=6t$, $x=6t-8$ and $x=-12t+9$ are the only $x$-coordinates of degree 1 defined over $\Q$. These sections are disjoint from the zero-section, and intersect at $t=\infty$ the singular fiber in the non-identity component. One can choose the $y$-coordinates of the three sections in such a way that the third sections is disjoint from the first two. This implies that the height pairing (see  \cite[Definition 8.5]{Sd}) yields the following intersection matrix
\[ \left( \begin{array}{ccc} {4}/{3}  & {2}/{3}& {2}/{3}\\ {2}/{3}& {4}/{3} &{2}/{3} \\ {2}/{3}&{2}/{3}& {4}/{3} \\\end{array} \right).\]
The determinant of this matrix is non-zero. This implies that these three sections generate a rank 3 subgroup $G_1$, and $G_1=MW_{\Q}(\pi)$.

Multiplying the $x$-coordinates with a cube root of unity will yield another rank 3 subgroup $G_2$ corresponding to a different eigenspace of the action of complex multiplication hence  $G_1\oplus G_2$ generate a rank 6 subgroup of the Mordell-Weil group.
\end{proof}
\section{Explicit formulas}\label{results}
Fix a field $K$ of characteristic 0. Let $E: y^2=x^3+ax+b$ and $F:y^2=x^3+cx+d$ be elliptic curves, with complete two-torsion over $K$. Assume that $E$ and $F$ have distinct $j$-invariant. Fix $\lambda$, $\mu$, $\nu$ and $\xi$ such that $E$ is isomorphic to $y^2=x(x-\lambda)(x-\mu)$ and $F$ is isomorphic to $y^2=x(x-\nu)(x-\xi)$. Then we may assume that 
\[ a=\frac{1}{3}(\lambda\mu-\lambda^2- \mu^2),\; b=\frac{1}{27} (3\lambda \mu(\lambda+\mu)-2(\lambda^3+\mu^3)),\]
and similar equations for $c$ and $d$.

The Kummer surface $Y$ is birational
to the surface $S\subset \A^3$ given by $(x^3+ax+b)t^2=z^3+cz+d$, and the fibration
$\psi:Y\ra \Ps^1$ is corresponds to the map $(x,z,t)\mapsto t$.

\begin{Lem}\label{LemTF} The groups $MW(\pi_i), i\geq 1$ are torsion-free.\end{Lem}
\begin{proof} Kuwata \cite[Theorem 4.1]{Kuw} shows that $\pi_{6}$ has smooth fibers over $t=0,\infty$ and  only singular fibers of type $I_1$ or $II$.
Hence the same holds for $\pi_{6i}$. This fact together with  \cite[Corollary VII.3.1]{MiES} implies that the group $MW(\pi_{6i})$ is
  torsion-free. Since $MW(\pi_{i})$ is a subgroup of $MW(\pi_{6i})$ it is
  also torsion-free.\end{proof}

We now discuss how to find explicit formulas for generators of $MW(\pi_i), i=1,\dots 6$. 
\subsection{$\pi_1: X_1\ra\Ps^1$}\label{ssecone}
It is not easy to find explicit equations for sections on
$\pi_1$, since there are infinitely many cases, depending on the minimal degree of an isogeny $E\ra F$. Instead we give a sufficient condition to have rank 1 or 2 over $K$.

\begin{Lem} The group $MW_K(\pi_1)$ has rank at
  most 2. 
\begin{itemize}
\item If $E$ and $F$ are isogenous over $K$ then $MW_K(\pi_1)$ has  positive rank. 
\item If $E$ and $F$ are isogenous over $K$ and $E$ admit complex multiplication then \\ $MW_K(\pi_1)$ has rank 2. 
\item If $MW_K(\pi_1)$ has positive rank then there exists a degree at most two extension $L/K$ such that $E$ and $F$ are isogenous.
\item if $MW_K(\pi_1)$ has rank 2 then there exists a degree at most two extension $L'/L$ such that $E$ admits complex multiplication over $L'$.
\end{itemize}
\end{Lem}

\begin{proof} Since $E$ and $F$ have complete two torsion it follows that $\rank MW_K(\pi_1)=r$ if and only if $18+r=\rank NS_K(X_2)=\rank NS_K(Y)$, with $Y$ the Kummer surface of $E$ and $F$.

It is easy to see that if $E$ and $F$ satisfy the first, resp., second assumption then the rank of $ NS_K(Y)$ is at least 19, resp., 20.

If $NS_K(Y)\geq 19$ then  $E$ and $F$ are isogenous over some extension of $K$. Let $\Gamma$ be the graph of the isogeny. Then the push-forward of $\Gamma$ on $Y$ is Galois-invariant. This implies that $([-1]\times [-1])^*\Gamma+\Gamma$ is Galois-invariant. From this it follows that $E$ and $F$ are isogenous over a degree 2 extension.

If $NS_K(Y)=20$ then $E$ and $F$ are isogenous over some extension of $K$ and $E$ has potential complex multiplication. Let $\Gamma$ be the graph of an isogeny, let $\Gamma'$ be the graph of the isogeny composed with complex multiplication. Then the push-forward of both $\Gamma$ and $\Gamma'$ are Galois-invariant. As above, this implies that $E$ admits complex multiplication over a degree 2 extension of $L$.
\end{proof}

\begin{Rem}\label{notisogexa} Suppose $E$ and $F$ are elliptic curves, not isogenous over $K$, but isogenous over a degree 2 extension $L/K$. Suppose that $\End(E)=\Z$. Let $\varphi: E \ra F$ be an isogeny defined over $L$. It is an easy exercise to show that the divisor 
\[ D:= \{ (P,\varphi(P)) | P \in E\} \cup \{ (P,-\varphi(P)) \mid P\in E\}\subset E \times F \]
is invariant under the action of $\Gal(L/K)$. Hence the push-forward of $D$ onto the Kummer is invariant under the action  of $\Gal(\bar{K}/K)$.  The argument used in the above proof  gives that $MW_K(\pi_1)$ has rank at least 1, while $E$ and $F$ are not isogenous over $K$.
\end{Rem}


\subsection{$\pi_2: X_2\ra\Ps^1$}\label{ssectwo} Since $MW(\pi_2)$ is torsion-free, we have that $MW(\pi_1') \oplus MW(\pi_1)$ is of finite index in $MW(\pi_2)$, where $\pi_1':X_1'\ra \Ps^1$ is the twist of $\pi_1$ at $0$ and $\infty$. 

One can easily show that $MW(\pi_1')\subset MW(\pi_3')$, where $\pi_3'$ is the twist of $\pi_3$ at $0$ an $\infty$. Since $MW(\pi_3')$ can be considered in a natural way as a subgroup of $ MW(\pi_6)$, we refer to that subsection for a discussion of the results. The explicit equations for generators of $MW(\pi_2')$ are given in \cite{Worksheet}. If we drop for a moment the condition  that $E$ and $F$  have complete 2-torsion over $K$, the we can use the result in \cite{Worksheet} to prove  that $MW_{K}(\pi_1')$ has rank 4 if and only if $E$ and $F$ have complete two-torsion over $K$.

\subsection{$\pi_3: X_3\ra\Ps^1$}\label{ssecthree}
A Weierstrass equation for $\pi_3$ is
\[y^2=x^3 -48acx +( \Delta(F)t^3+864bd+\Delta(E)t^{-3}).\]
Setting $s=(t+\alpha_i/t)$, with $\alpha_i^3=\Delta(E)/\Delta(F)$ and $i=1,2,3$, will give an equation for the rational elliptic surface $\psi_i: S_i\ra \Ps^1$, given by
\[y^2=x^3-48acx - \Delta(F) (s^3-3\alpha_i s) +864bd.\]
This surface has a fiber of type $I_0^*$ at $s=\infty$. 
For each choice of $\alpha_i$ we obtain an isomorphic surface over $K(\zeta_3)$, but the pullback of the sections to $\pi_3$ depends on the choice of $\alpha_i$: 

\begin{Lem}\label{LemOrt} For $i\neq j$ we have $MW(\psi_i) \cap MW(\psi_j)=\{\sigma_0\}$, where $MW(\psi_i)$ and $MW(\psi_j)$ are considered as subgroups of $MW(\pi_3)$.\end{Lem}

\begin{proof} Without loss of generality, we may assume that $i=1$ and $j=2$. Set $L:=K(\zeta_3,\alpha_1)$. Let $G\subset\Aut(L(t))$ be generated by $\varphi_1:t\mapsto \alpha_1/t$ and $\varphi_2:t \mapsto \alpha_1\zeta_3/t$.
From $\varphi_1\circ \varphi_2: t \mapsto t\zeta_3$ it follow that $\#G\geq 6$. Set $s'=t^3+\alpha_1^3/t^3$. Then $s'$ is fixed under $G$. We have the following inequalities
 \[6\leq \#G =[L(t)^G:L(t)] \leq [L(s'):L(t)]=6.\] These inequalities give that $L(t)^G=L(s')$. This implies that a section in $MW(\psi_i)\cap MW(\psi_j)$ is the pull back of a section of the elliptic surface $\psi'$ with Weierstrass equation
\[y^2=x^3-48acx-\Delta(F) s'+864bd.\]
This is an equation of rational elliptic surface with a $II^*$-fiber. In this case the Shioda-Tate formula~\ref{ST} implies that $MW(\psi')$ has rank 0.  Since $MW(\psi')$ is a subgroup of the torsion-free group $MW(\pi_3)$ (see Lemma~\ref{LemTF}) it follows that $\#MW(\psi')=1$.
\end{proof}

\begin{Lem}\label{LemFour} We have that $\rank MW(\psi_i)=4$.\end{Lem}
\begin{proof} From the equation of $\psi_i$ one easily sees that it a rational elliptic surface, with a fiber of type $I_0^*$ over $s=\infty$, and no other reducible singular fibers. Hence the Shioda-Tate formula~\ref{ST} implies that $\rank MW(\psi_i)=4$.\end{proof}

\begin{Lem}\label{downstairs} For $i\neq j$ we have $(MW(\psi_i)\oplus MW(\psi_j))\cap MW(\pi_1)=\{ \sigma_0\}$, considered as subgroups of $MW(\pi_3)$.\end{Lem}

\begin{proof} From Lemma~\ref{LemOrt} and Lemma~\ref{LemFour} it follows that $MW(\psi_i)\oplus MW(\psi_j)$ injects into $MW(\pi_3)$. Consider the vector spaces $V=(MW(\psi_i)\oplus MW(\psi_j))\otimes \C$ and $W=MW(\pi_1)\otimes \C$. The automorphism $\sigma:(x,y,t)\mapsto (x,y,\zeta_3t)$  induces a  trivial action on $W$. We now prove  that $\sigma$ maps $V$ to itself and all eigenvalues of this action are different from 1.

Assume for the moment that $\lambda,\mu,\nu$ and $\xi$ are algebraically independent over $\Q$. This defines a  Kuwata elliptic surface
\[ \Pi_k: \mathcal{X}_k \ra \Ps^1_{\Q(\lambda,\mu,\nu,\xi,\sqrt[3]{\lambda(\lambda-\mu)\mu/\nu(\nu-\xi)\xi})} \]
and, similarly,
\[ \Psi_k: \mathcal{S}_k \ra \Ps^1_{\Q(\lambda,\mu,\nu,\xi,\sqrt[3]{\lambda(\lambda-\mu)\mu/\nu(\nu-\xi)\xi})}.\]

Since $MW(\Pi_1)=0$, and $MW(\Psi_i)\oplus MW(\Psi_j)$ has rank 8, it follows that a lift $\tilde{\sigma}$ of $\sigma$ acts on $\tilde{V}:=MW(\Psi_i) \oplus MW(\Psi_j)$, hence $\sigma$ acts on $V$. If $\tilde{\sigma}$ would have an eigenvalue 1 on $\tilde{V}$, then $MW(\Pi_1)$ would be non-trivial. Hence also $\sigma$ acts without eigenvalue $1$. This implies that $V\cap W=0$, yielding the lemma.
 \end{proof}
From the Lemmas~\ref{LemOrt},~\ref{LemFour} and~\ref{downstairs} and Theorem~\ref{exactrank} it follows that for $i\neq j$ we have that $MW(\psi_i)\oplus MW(\psi_j)$ generates a subgroup of finite index in $MW(\pi_3)/MW(\pi_1)$. Hence if one wants  to describe the Galois representation on $MW(\pi_3)/MW(\pi_1)$  it suffices to describe the Galois representation on $MW(\psi_i)$.

Since $\psi_i:S_i \ra \Ps^1$ are rational elliptic surfaces with an additive fiber at $s=\infty$ we can apply Section~\ref{secfind} to calculate expressions for the generators. The relevant formulae for this may be found in \cite{Worksheet}.

\subsection{$\pi_4: X_4\ra\Ps^1$}\label{ssecfour}

Let $\pi'_2:X'_2\ra \Ps^1$ be the twist of $\pi_2$ by the points $0$ and $\infty$. 
 The results of Section~\ref{Prel} imply that
the group $MW(\pi_2)\oplus MW(\pi_2')$ is of finite index in $MW(\pi_4)$. Since the fibers over $0$ and $\infty$ of $\pi_2$ are of
type $IV^*$ and all other fibers are of type $I_1$ or $II$, it follows from the results mentioned in Section~\ref{Prel} that
$\pi'_2$ has only fibers of type $II$ or $I_1$. Moreover, $\pi_2':X_2'\ra \Ps^1$ defines a rational elliptic surface, hence $\rank MW(\pi_2')=8$.

\begin{Lem}\label{rang2} Let $\pi: X \ra \Ps^1$ be a Jacobian elliptic surface 
  such that all singular fibers are irreducible. Let
  $Z:=\sigma_0(\Ps^1)$. Let $S$ be the image of a section
  $\sigma:\Ps^1\ra X$. Denote $S_n$ the image of $(n\sigma): \Ps^1\ra
  X$.

Then the equality $(S_n\cdot Z)=(S_m\cdot Z)$ holds if and only if $n=\pm m$.\end{Lem}

\begin{proof} Let $\langle\cdot ,\cdot \rangle$ denote the height paring on $MW(\pi)$ (see
  \cite[Definition 8.5]{Sd}). In this case $\langle T,T \rangle=2(T\cdot Z)-2\chi(X)$, hence
\[ 2(S_n\cdot Z)-2\chi(X) =\langle S_n,S_n \rangle = \langle nS_1,nS_1
\rangle =n^2 \langle S_1,S_1 \rangle =2n^2(S_1\cdot Z)-2n^2\chi(X).\]
It follows that $(S_n\cdot Z)=n^2((S\cdot  Z)-\chi(X))+\chi(X)$, which yields the lemma.
\end{proof}

Assume, for the moment, that $\lambda,\mu,\nu$ and $\xi$ are algebraically independent over $\Q(t)$.
We can consider the elliptic surface  $\pi_2'$ as an elliptic curve $A$ over $K':=\Q(\lambda,\mu,\nu,\xi,t)$.

Then on $A$  we have that
\begin{small} { \[ x= \frac{2(\sqrt{\lambda}+\sqrt{\mu})\lambda\mu}{\sqrt{\xi}+\sqrt{\xi-\nu}} \frac{1}{t} +\left(-\frac{4}{3} (2\xi-\nu)(\lambda+\mu)+4\sqrt{\lambda\mu\xi(\xi-\nu)}\right) + \frac{2(\sqrt{\xi}+\sqrt{\xi-\nu})\xi (\xi-\nu)}{\sqrt{\lambda}+\sqrt{\mu}}t\]}
\end{small}
\noindent is a $x$-coordinate of a point $P_1$, hence giving rise to 2 different
points on $A$ (see \cite{Worksheet}). If we plug this $x$-coordinate in the equation for $A$, then
\[ y^2= 2(\sqrt{\xi}+\sqrt{\xi-\nu})(\sqrt{\lambda}+\sqrt{\mu}) h(t)^2 \]
for some $h\in \Q(t,\sqrt{\lambda},\sqrt{\mu},\sqrt{\xi},\sqrt{\nu-\xi})$ (see \cite{Worksheet}). 

Using $P_1$ we now show how to find points $P_2,\dots P_8$, such that if we specialize, the $P_i$ generate $MW(\pi_2')$. 
Let $P_2$ be a section such that $x(P_2)$ is obtained from $x(P_1)$ by replacing $\sqrt{\mu}$ by $-\sqrt{\mu}$. (The existence follows from the fact that the map $\sqrt{\mu}\mapsto -\sqrt{\mu}$ fixes the equation of $\pi_2'$.)

\begin{Lem} The points $P_1$ and $P_2$ are independent in $A(\bar{K})$.
\end{Lem}

\begin{proof} Since $x(P_1)\neq x(P_2)$ it follows that $P_1\neq \pm P_2$. Hence the sections $S_1$ and $S_2$ on $\pi_2'$ satisfy  $(S_1\cdot Z)=(S_2\cdot Z)$ and $S_1\neq \pm S_2$. It follows 
from Lemma~\ref{rang2} that $S_1$ and $S_2$  generate a subgroup of $MW(\pi'_2)$ of rank 2, hence $P_1$ and $P_2$ are independent.
\end{proof}


Note that we have $P_1\in A(K_+)$ and $P_2\in A(K_-)$, where 
\[ K_{\pm}:=K'(\sqrt{\lambda},\sqrt{\mu},\sqrt{\xi},\sqrt{\xi-\nu},\sqrt{ 2(\sqrt{\xi}+\sqrt{\xi-\nu})(\sqrt{\mu}\pm\sqrt{\lambda})} ).\]

The automorphism $\sigma$ of $K'$, fixing $\Q$ and further given by \[(t,\lambda,\mu,\nu,\xi)\mapsto \left(\frac{1}{t},\nu,\xi,\lambda,\mu\right)\]
is an automorphism of $A$, (i.e., swapping $E$ and $F$, and applying the $\Ps^1$-automorphism $t\mapsto 1/t$). This yields to different points $P_3:=\sigma(P_1),P_4:=\sigma(P_2)$ on $A$.

\begin{Lem} The points $P_1,\dots P_4$ generate a rank 4 subgroup of $A(\bar{K})$.\end{Lem}

\begin{proof} The points $P_3$ and $P_4$ are not in $A(K_+K_-)$, but are contained in $A(K'_{\pm})$, where \[K'_{\pm}= K'(\sqrt{\mu},\sqrt{\mu-\lambda},\sqrt{\nu},\sqrt{\xi},\sqrt{ 2(\sqrt{\mu}+\sqrt{\mu-\lambda})(\sqrt{\xi}\pm\sqrt{\nu})} ).\]
Applying \cite[Lemma 1.3.2]{Scholten} yields that the $P_i$ generate a rank 4 subgroup.
\end{proof}

Let $\tau'$ automorphism of $K'$ fixing $\Q(t)$ and mapping $(\xi,\lambda)\leftrightarrow(\nu,\mu)$ (i.e., we are interchanging the role of the two-torsion points of both $E$ and $F$). Set $P_i'=\tau(P_i)$.

\begin{Lem} The points $P_1,\dots, P_4,P'_1,\dots,P'_4$ generate a subgroup of rank 8 of $A(\bar{K'})$.
\end{Lem}

\begin{proof}
The points $P_1',\dots, P_4'$ are defined over different fields (i.e., they are defined over $K_{\pm}(i)$ or $K_{\pm}'(i)$, but not over $K_{\pm}$ or $K_{\pm}'$.) From \cite[Lemma 1.3.2]{Scholten} it follows then that the rank of the subgroup generated by the $S_i, S'_i$ is 8.
\end{proof}
 
\begin{Prop} The sections of $\pi_2'$ associated to the points $P_1,\dots P_4,P_1',\dots P_4'$ of $A$ generate a rank 8 subgroup of $MW(\pi_2')$.\end{Prop}

\begin{proof} Let $\Pi'_2:\mathcal{X}'_2 \ra \Ps^1_{\Q(\lambda,\mu,\nu,\xi)}$ be the elliptic surface corresponding to $A/K'$. The points $P_i,P'_i$ define sections $\mathcal{S}_i,\mathcal{S}'_i$ of $\Pi_2'$, which generate a rank 8 subgroup of $MW(\Pi_2')$. This implies that for a general specialization the specialized sections $S_i, S_i'$ generate a rank 8 subgroup of $MW(\pi_2')$. Since the intersection numbers are constant under deformation, it follows that if the $S_i, S_i'$ do not generate $MW(\pi_2')$, then $\pi_2'$ has at least one reducible fiber. From \cite[Theorem 4.3]{Kuw} it follows that then $E$ and $F$ are isomorphic, which contradicts our assumption at the beginning of this section.
\end{proof}

\begin{Cor}\label{res4}  The sections of $\pi_2'$ associated to the points $P_1,\dots,P_4,P_1',\dots,P_4'$ generate a subgroup of finite index in $MW(\pi_4)/MW(\pi_2)$.\end{Cor}

\begin{proof} This follows from the above Proposition together with the observation that the direct sum of  $MW(\pi_2')$ and $ MW(\pi_2)$ generates a subgroup of finite index of $MW(\pi_4)$.
\end{proof}

We will now take special values for $\lambda,\mu,\nu,\xi$ such that several sections are defined over $\Q$. 

\begin{Cor}\label{CorQ} For infinitely many pairs of elliptic curves $(E,F)$ over $\Q$, the group $MW_{\Q}(\pi_4)$ has rank 8 or 9.\end{Cor}

\begin{proof} We have that $\rank MW_{\Q}(\pi_4)=\rank MW_{\Q}(\pi_2)+\rank MW_{\Q}(\pi_2')$. From the results in \ref{ssecone} and \ref{ssectwo} we know that $\rank MW_{\Q}(\pi_2)$ is either 4 or 5.

As mentioned in Section~\ref{notres}, one can prove that the rank over $\Q$ of $\pi_2'$ is at most 4. 
To obtain rank 4 it suffices to choose $\lambda,\mu,\nu,\xi$ such that $S_1, S_2, S_3$ and $S_4$ are defined over $\Q$. In order to obtain this set 
 $\lambda=l^2,\mu=m^2,\nu=n^2$, $\xi=k^2$, $k^2-n^2=n_2^2$ and $m^2-l^2=l_2^2$. 

Then
\[ 2(k+n_2)(m\pm l), 2(m+l_2)(k\pm n) \]
have to be a non-zero square for all choices of $\pm$. One easily computes that this occurs precisely when
\[l_2n_2=\frac{u^2(\rho-1)(\tau^2-1) }{4(\rho+1) \tau^2} \]
for some $u\in \Q^*$ and
\[
k= n_2\frac{ \tau^2+1}{\tau^2-1}, \;\;\;\;
n= n_2 \frac{2 \tau}{\tau^2-1},\;\;\;\;
m = l_2 \frac{\rho^2+1}{\rho^2-1},\;\;\;\;
l = l_2 \frac{2 \rho}{\rho^2-1}.\]
From these last equations we can obtain our original $\lambda,\mu,\nu,\xi$. 

The associated Legendre parameters are
\[\left( \frac{2 \tau}{\tau^2+1}\right)^2 \mbox{ and } \left(\frac{2 \rho}{\rho^2+1}\right)^2.\]
One easily finds $\rho,\tau$ such that the corresponding $j$-invariants are different.\end{proof}

For later use, we remark that $\Delta(E)/\Delta(F)$ is a third power if and only if 
\[\frac{\tau (\tau^4-1)}{\rho(\rho^4-1)}\] is a third power. 
A necessary condition for obtaining curves with different $j$-invariant is $\rho\neq \tau$. The only solution we found with these properties is $(\rho,\tau)\in \{(2,3),(3,2)\}$. This gives rise to two curves having both  Legendre parameter $25/9$, hence does not give an interesting solution.

\subsection{$\pi_5: X_5 \ra \Ps^1$.} \label{ssecfive}
An equation for $\pi_5$ is
\[y^2=x^3-48acx +( \Delta(F)t^5+864bd+\Delta(E)t^{-5})\]
Setting $s=(t+\alpha_i/t)$, with $\alpha_i^5=\Delta(E)/\Delta(F)$, gives  a rational elliptic surface $\psi_i: S_i \ra \Ps^1$.

We can now copy the strategy used for  of $\pi_3$. We have that $MW(\pi_5)/MW(\pi_1)$ has rank  16, the $MW(\psi_i)$ have rank 8, the intersection $MW(\psi_i)\cap MW(\psi_j)=\{\sigma_0\}$, considered as subgroups of $MW(\pi_5)$ and $(MW(\psi_i)\oplus MW(\psi_j))\cap MW(\pi_1) = \{\sigma_0\}$. This combined with the fact that all the $\psi_i$ are isomorphic over $K(\zeta_5)$ implies that  it suffices to find sections of only one of the $\psi_i$. One can show that $\psi_i$ is a rational elliptic surface with an additive fiber at $t=\infty$. Hence we can apply Section~\ref{secfind} to find an expression for the sections of $MW(\pi_5)$. In all choices  for  $\lambda,\mu,\nu,\xi$ we tried the final resultant is a product of two polynomials of degree 120. (See \cite{Worksheet})

\subsection{$\pi_6: X_6 \ra \Ps^1$.}\label{ssecsix}
Since $\pi_6$ is torsion-free we have that $MW(\pi_3)\oplus MW(\pi_3')$ is of finite index in $MW(\pi_6)$, where $\pi_3':X'_3 \ra \Ps^1$ is the twist of $\pi_3$ at 0 and $\infty$. 

The group $MW(\pi_3)$ is described above. We discuss here how to find generators for $MW(\pi_3')$. We follow Kuwata \cite[Section 5]{Kuw}.
 Kuwata observes  that $X'_3$ is birational over $\Q$ to the cubic surface in $\Ps^3$
\[ C:Z^3+cZY^2+dY^3=X^3+aXW+bW^3,\]
and that  the strict transforms of the 27 lines of $C$ to $X'_3$ generate $MW(\pi_3')$.

It has been known for a long time how to find the 27 lines on a cubic surface of the above form, see for example \cite{Schlaefli}.
 Due to our special situation  we give a somewhat different approach to find all 27 lines.

In our case $C$ is isomorphic to
\[ Z(Z-\nu Y)(Z-\xi Y)=X(X-\lambda W)(X-\mu W).\tag{*}\]
One finds  9 lines defined over $\Q$, namely the intersections of $Z-\alpha Y=0$ and $X-\beta W=0$, with $\alpha \in \{0,\nu,\xi\}$ and $\beta \in \{0,\lambda,\mu\}$. We give now the equations for the strict transforms of these lines on $X_3$: 
\begin{Lem} Set
\[ A =(a_{ki})= \left( \begin{array}{ccc}
\nu-\xi&\xi-\nu&\nu\\
-\xi&-\nu&\xi\\
-\mu& -\lambda& \lambda\\
\lambda-\mu& \mu-\lambda&\mu\end{array}\right).
 \]
Then 
\[x=4a_{1i}a_{2i} t^2+\frac{4}{3} (a_{1i}+a_{2i}+a_{3j}+a_{4j})t+4a_{3j}a_{4j}\] and \begin{small}
\[y=4a_{1i}a_{2i}(a_{1i}+a_{2i})t^3+8a_{1i}a_{2i}(a_{3j}+a_{4j})t^2+8a_{3j}a_{4j}(a_{1i}+a_{2i})t+4a_{3j}a_{4j}(a_{3j}+a_{4j})\]\end{small}
for $1\leq i,j\leq 3$ are sections of $\pi'_3$. 
\end{Lem}

\begin{proof} This is a straightforward computation.
See \cite{Worksheet}.
\end{proof}

The usual strategy to find the other 18 lines is to take pencils of planes through the lines we found above. We choose a different strategy. The image of a section on $X_3'$ can be pulled back to a divisor on  $X_6$. We can take the push-forward to $Y$ of this divisor and then push-forward in a natural way to $\Ps^1\times \Ps^1$. We give 6 divisors on $\Ps^1\times \Ps^1$, such that the pull back of such a divisor to $X_6$ consists of three components, and the 18 divisors obtained in this way correspond to the 18 lines on $C$.
Consider $\Ps^1 \times \Ps^1$ with projective coordinates $X,W$ and $Y,Z$. Set 
\[ P_1:=[0,1], P_2:= [\lambda,1], P_3:=[\mu,1], Q_1:=[0,1], Q_2:=[\nu,1], Q_3:=[\xi,1].\]
Fix a permutation $\sigma \in  S_3$. Let $C_{\sigma} \subset \Ps^1\times \Ps^1$ be the $(1,1)$-curve going through $(P_i,Q_{\sigma(i)})$. For example, $C_{\id}$ is given by
\[G(X,Y,Z,W):= \nu\xi(\lambda-\mu)X W +\mu\lambda(\xi-\nu)WZ+ (\nu\mu -\xi\lambda)XZ =0.\]

We describe what the corresponding divisor on  $C$ is. This can be done
by dehomogenizing first, say, we set $W=Y=1$. Then we can express $Z$  in terms of $X$. Substitute this expression for $Z$  in $(*)$. We obtain a rational expression in $X$ (a quotient of a polynomial of degree 6 and a polynomial of degree 3). The numerator is zero if and only if $X=0,X=\lambda, X=\mu$ or $X$ satisfies a degree 3 polynomial $f$. For example if $\sigma=\id$ then 
\[\left(\frac{\mu\lambda(\nu-\xi) + \sqrt[3]{\mu\lambda (\mu-\lambda) \nu^2\xi^2(\nu-\xi)^2}}{\nu\mu-\xi\lambda} \right)\zeta_3^k=:\gamma_k
\]
are the zeroes of $f$. Let $H$ be the hyperplane $X=\gamma_k Y$. Let $H'$ be $G(X,Y,\gamma_k, 1)=0$. Then $H\cap H'$ is contained in $C$ and defines a line $\ell_{\sigma,k}$.
 One can easily show that for all $\sigma \in H$ one obtains three such hyperplanes, and that one obtains 18 lines in total. 

\begin{Rem}One can  show that 6 of these 18 lines are defined over 
\[ K(E[2],F[2],\sqrt[3]{\Delta(E)/\Delta(F)}),\]
and the 12 others over 
\[ K(E[2],F[2],\mu_3,\sqrt[3]{\Delta(E)/\Delta(F)}).\] Moreover, one shows easily that if one of the 18 lines is defined over $K(E[2],F[2],\mu_3)$, then the quotient $\Delta(E)/\Delta(F)$ is a third power $K(E[2],F[2])$. This contradicts the claim in \cite[Section 5]{Kuw}, which states that all  27 lines are defined over $K(E[2],F[2],\mu_3)$.\end{Rem}

\begin{Rem} One can easily find elliptic curves  $E,F$ over $\Q$ with complete two-torsion defined over $\Q$ such that $\Delta(E)/\Delta(F) \in \Q^{*3}$: one needs to find solutions of
\[ \lambda(\lambda-\mu)\mu=\tau^3 \nu(\nu-\xi)\xi. \]
This defines a conic over $\Q(\tau,\mu,\xi)$ containing the rational point $\lambda=0$ and $\nu=0$. Hence one can parameterize this conic and find solution giving rise to elliptic curves. For example $(\lambda,\mu,\nu,\xi,\tau)=(16,1,6,1,2)$ gives an example.
\end{Rem}

\begin{Rem} Suppose one can find a  solution in $\Q^3$ of
\[ \frac{\tau(\tau^4-1)}{\rho(\rho^4-1)}= \sigma^3\]
satisfying the conditions mentioned at the end of \ref{ssecfour}. Then we obtain examples such that  $\rank MW_{\Q}(\pi_6)$ is at least $ 6$ and $\rank MW_{\Q}(\pi_4)$ is at least 8, yielding  that the rank of $MW_{\Q}(\pi_{12})$ is at least $10$.
\end{Rem}


\end{document}